\documentclass[preprint,12pt]{elsarticle}

\usepackage{amssymb}
\usepackage{amsmath}
\usepackage{amsmath}
\usepackage{amsthm}
\usepackage{amssymb}
\usepackage{tikz-cd}
\newtheorem{Theorem}{Theorem}[section]
\newtheorem{lemma}[Theorem]{Lemma}

\newtheorem{definition}[Theorem]{Definition}
\newtheorem{corollary}[Theorem]{Corollary}
\newtheorem{proposition}[Theorem]{Proposition}
\newtheorem{example}[Theorem]{Example}
\newtheorem{remark}[Theorem]{Remark}

\journal{arXiv}

\begin{document}

\begin{frontmatter}


\author{Reza Ameri\fnref{label1}}
\ead{ rameri@ut.ac.ir}

\author{Behnam Afshar\fnref{label1}}
\ead{behnamafshar@ut.ac.ir}

\affiliation[label1]{organization={Department of Mathematics, Statistics and Computer Science, College of Science, University of Tehran},
	            city={Tehran},
                postcode={14155-6455},
	            country={Iran}}

\title{Zariski topology of (Krasner) hyperrings}




\begin{abstract}
In this article, we will study prime spectrum of Krasner hyperrings and Zariski topology on them, which play an important role in algebraic geometry. Then some results about the relationship between the topological properties of $Spec(R)$ and the algebraic properties of the hyperring $R$ will be proved. In the following, by proving that every strongly regular relation on Krasner hyperrings can be considered as a congruence relation, we will define a topology on the set of strongly regular relations, and investigate its relationship with the Zariski topology. In addition, the effect of fundamental relations on the Zariski topology of Krasner hyperrings will also be investigated.
\end{abstract}




\end{frontmatter}




\section{Introduction}    \label{sec1}
The theory of algebraic hyperstructures was first proposed in 1934 by F. Marty, who defined hypergroups as a generalization of groups \cite{bib17}. F. Marty showed that the characteristics of hypergroup can be used in group theory, rational functions and algebraic functions. Surveys of the theory can be found in \cite{bib10,bib11}. Later, hyperstructure theory was studied in connection with various fields. A. Connes and C. Consani \cite{bib6,bib7,bib8} have shown that there is a relationship between the hyperstructure theory and algebraic geometry.\\
In 1956 M. Krasner introduced the concept of hyperrings to use them as a technical tool on the approximation of valued fields \cite{bib16}. A General hyperring is a hyperstructure $(R,+,.)$ where $(R,+)$ is a hypergroup and $(R,.)$ is semihypergroup such that ' . ' is distributive with respect to '+'. If $(R,+)$ is a canonical hypergroup and $(R,.)$ is semigroup such that zero element is absorbing, then $(R,+,.)$ is called the Krasner hyperring \cite{bib9,bib10}.\\
The fundamental relations are
one of the most important and interesting concepts in algebraic hyperstructures
that ordinary algebraic structures are derived from algebraic hyperstructures by
them. The fundamental relation $\beta^{\ast}$ on hypergroups was defined by Koskas \cite{bib15},
Corsini \cite{bib9}, Ferni \cite{bib14}, and Vogiouklis \cite{bib21}. Then D. Ferni introduced the
fundamental relation $\gamma^{\ast}$ which is the transitive closure of $\gamma$ and is the smallest
relation such that $H/\gamma^{\ast}$ is an abelian group.
T. Vogiouklis generalized the fundamental relations in \cite{bib21} to use on hyperrings
and in \cite{bib1}, it has been demonstrated that relations $\beta^{\ast}$ and $\gamma^{\ast}$ are related together
in the form of $\gamma^{\ast}=\delta\ast\beta^{\ast}$, where $\delta$ is the congruence relation with respect to the
commutator subgroup.\\
Also, in \cite{bib2}, a one-to-one correspondence between the lattice of strongly regular relations on regular hypergroup $H$, and the lattice of normal subhypergroups of $H$ containing $\omega_{H}$, is introduced and in \cite{bib4} the fundamental relation on hypermodules is introduced.\\
Recently, many studies have been conducted on the properties of Zariski topology of algebraic hyperstructures. To study the Zariski topology of multiplicative hyperrings, refer to \cite{bib3}. Here we are going to investigate the characteristics of Zariski topology of Krasner hyperrings, which is used in the study of sheaves of hyperrings.
The main features of Zariski topology of Krasner hyperring $R$, will be discussed in section \ref{sec3}, and we will show that the necessary and sufficient condition for $Spec(R)$ to be connected is that $R$ is not nontrivial product of hyperrings. We also show that the necessary and sufficient condition for $Spec(R)$ to be irreducible is that $nil(R)$ is prime hyperideal of $R$. In the following, other properties such as compactness and chain conditions will be examined for $Spec(R)$. 
In section \ref{sec4}, we define the Zariski topology functor $Spec(-)$ for hyperrings and prove some of its categorical properties. Also with the help of equivalence relation $\gamma^{*}$ we will examine the relationship between this functor and its classical Zariski topology functor.
Furthermore, a one-to-one correspondence between strongly regular relations and hyperideals containig $\gamma^{\ast}(0)$ will be introduced. Also, by introducing prime and primary strongly regular relations, a topology is defined on the set of all strongly regular relations, which in the next researches we intend to examine its connection with Zariski topology.

\section{Preliminaries}     \label{sec2}
A Krasner hyperring is an algebraic structure $(R,+,.)$ which satisfies the following axioms:
\begin{itemize}
	\item{For every $x,y,z\in R$, $x+(y+z)=(x+y)+z$;}
	\item{For every $x,y\in R$, $x+y=y+x$;}
	\item{There exists $0\in R$ such that $0+x=\{x\}$, for every $x\in R$;}
	\item{For every $x\in R$ there is a unique $x'\in R$ that $0\in x+x'$(we use $-x$ for $x'$);}
	\item{If $z\in x+y$ then $y\in-x+z$ and $x\in z-y$.}
\end{itemize}
A Krasner hyperring $(R,+,.)$ is called commutative (with unite element) if $(R,.)$ is commutative (with unite element) semigroup. If $(R-\{0\},.)$ is a group then $(R,+,.)$ is called a Krasner hyperfield and if $(R,+,.)$ is commutative Krasner hyperring with unite element and $ab=0$ implies that $a=0$ or $b=0$ for all $a,b\in R$, then $R$ is called hyperdomain. Throughout this paper, hyperring refers to commutative Krasner hyperring with a unit element.

\begin{definition}[\cite{bib11}]           \label{ii-1}
	A subhyperring $I$ of $R$ is a left(right) hyperideal of $R$, If $ra\in I$ ($ar\in I$) for all $r\in R$, $a\in I$. If $I$ is both a left and a right hyperideal Then $I$ is called the hyperideal. A proper hyperideal $M$ of $R$ is maximal hyperideal of $R$, if the only hyperideals of $R$ contain $M$, are $M$ itself and $R$. Let $P$ be a proper hyperideal of $R$, then $P$ is called the prime hyperideal if for every pair of hyperideals $A$ and $B$ of $R$; $AB\subseteq P$ implies $A\subseteq P$ or $B\subseteq P$.
\end{definition}

\begin{lemma}[\cite{bib11}]     \label{ii-2}
	A nonempty subset $I$ of a hyperring $R$ is a left(right) hyperideal if and only if $a-b\subseteq I$ and $ra\in I$($ra\in I$), for all $a,b\in I$, $r\in R$. 
\end{lemma}

An equivalence relation $\theta$  on a Krasner hyperring $R$ is called \textit{regular} if the following implication holds:
\begin{equation}
	a\ \theta \ b \ ,\ c\ \theta \ d \ \Rightarrow (a+c)\ \bar{\theta}\ (b+d) \ \ and \ \ ac\ \theta \ bd,
\end{equation}
and is called \textit{strongly regular} if
\begin{equation}
	a\ \theta \ b \ ,\ c\ \theta \ d \ \Rightarrow (a+c)\ \bar{\bar{\theta}}\ (b+d) \ \ and \ \ ac\ \theta \ bd.
\end{equation}

\begin{definition}[\cite{bib11}]   \label{ii-3}
	Let $R$ and $S$ be hyperrings. A mapping $f$ from $R$ to $S$ is said to be a good homomorphism if for every $a,b\in R$: $f(a+b)=f(a)+f(b)$, $f(ab)=f(a)f(b)$ and $f(0)=0$. 
\end{definition}

Let $R$ be a commutative hyperring with a unit element and $I$ be a proper hyperideal of $R$. Then there exists a maximal hyperideal of $R$ containing $I$ and each maximal hyperideal is a prime hyperideal. Additionally, $I$ is prime if $ab\in P$ implies that $a\in P$ or $b\in P $, for every $a,b\in R$.

\begin{proposition}[\cite{bib11}]     \label{ii-6}
	Let $R$ be a commutative hyperring with a unit element and $I$ be a proper  hyperideal of $R$. Then:
		\item[(i)] $I$ is prime hyperideal if and only if $R/I$ is a hyperdomain;
		\item[(ii)] $I$ is maximal hyperideal if and only if $R/I$ is a hyperfield.
\end{proposition}

For any regular hypergroup $H$, if $\mathcal{SR}(H)$ is the set of all strongly regular relations on $H$ and $N(S_{\beta})$ is the set of all normal subhypergroups of $H$, containing $S_{\beta}$ ($=\omega_{H}$), then the map
\begin{equation}
	\underset{\rho \ \mapsto \ S_{\rho}}{\varphi: \mathcal{SR}(H)\rightarrow N(S_{\beta})}		
\end{equation}
where $S_{\rho}=\{x\in H;\ \rho(x)=e_{H/\rho}\}$, is an isomorphism of complete lattices \cite{bib2}.

\begin{definition}[\cite{bib19}]    \label{ii-12}
	The topological space $T$ is called a disconnected space if it can be decomposed as a disjoint union of two nonempty closed subsets, and is called the irreducible space if any pair of nonempty open subsets of $T$, intersect. A topological space $T$ is $T_{0}$-space if for any $a,b\in T$($a\neq b$), there exist open subsets $U$ and $V$ of $T$ such that $a\in U$, $b\notin U$ or $b\in V$, $a\notin V$. $T$ is $T_{1}$-space if for any $a,b\in T$($a\neq b$), there exists open subsets $U$ and $V$ of $T$ such that $a\in U$, $b\notin U$ and $b\in V$, $a\notin V$. $T$ is $T_{2}$-space(Hausdorff) if for any $a,b\in T$($a\neq b$), there exists open subsets $U$ and $V$ of $T$ such that $a\in U$, $b\in V$ and $U\cap V=\emptyset$.
\end{definition}

\begin{Theorem}[\cite{bib19}]     \label{ii-13}
	Let $T$ be a topological space. Then:
		\item[(i)] {If $S$ is an irreducible subspace of $T$ then $\bar{S}$ is irreducible;}
		\item[(ii)] {Every irreducible subspace is contained in a maximal irreducible subspace;}
		\item[(iii)] {The irreducible components of $T$ are closed and cover $T$.}
\end{Theorem}

\section{Zariski topology of Krasner hyperrings}\label{sec3}

Everywhere in this section $R$ is a commutative Krasner hyperring with a unit element.
Let $R$ be a hyperring and $Spec(R)$ be the set of all prime hyperideals of $R$ and $mSpec(R)$ be the set of all maximal hyperideals of $R$. For all $x=P\in Spec(R)$ let $\mathbb{K}(x)$ be the quotient hyperfield of hyperdomain $R/P$. For all $f\in R$ we have $R\rightarrow R/P\rightarrow \mathbb{K}(x)$ where $f\rightarrow f+P\rightarrow \frac{(f+P)(1+P)}{1+P}$.

\begin{remark}    \label{iii-2}
	Since every maximal hyperideal is a prime hyperideal, it is always true that $mSpec(R)\subseteq Spec(R)$. Also $R$ is hyperfield if and only if $0\in mSpec(R)$ if and only if the only hyperideals of $R$ are $0$ and $R$.
\end{remark}

\begin{definition}   \label{iii-3}
	For every subset $S\subseteq R$, let
	$$V(S)=\{x\in Spec(R) \ ; \ f(x)=0, \ \forall f\in S\}=\{P\in Spec(R)\ ; \ S\subseteq P\} $$
	$$V_{m}(S)=\{x\in mSpec(R) \ ; \ f(x)=0, \ \forall f\in S\}=\{M\in mSpec(R)\ ; \ S\subseteq M\}. $$
\end{definition}

\begin{lemma}    \label{iii-4}
	The following properties hold for $V$:
	\item[(i)] {Let $I=I(S)$, be the hyperideal generated by $S$, then $V(S)=V(I)$;}
	\item[(ii)] {If $S_{1}\subseteq S_{2}$, then $V(S_{2})\subseteq V(S_{1})$;}
	\item[(iii)] {$V(S)=\varnothing$ if and only if $1\in I(S)$;}
	\item[(iv)] {$V(M)=M$ if and only if $M\in mSpec(R)$.}
	\begin{proof}
		The proof is straightforward and similar to the classical case. 
	\end{proof}
\end{lemma}

\begin{Theorem}    \label{iii-5}
	The sets of family $\{V(I)\}_{I\lhd R}$, satisfy the axioms for closed sets in a topological space.
	\begin{proof}
		We have $V(R)=\emptyset$ and $V(\{0\})=Spec(R)$. Let $\{V(I_{j})\}_{j\in J}$ be a family of closed sets and $P\in V(\sum_{j\in J}I_{j})$. Then $I_{i}\subseteq \sum_{j\in J}I_{j}\subseteq P$, for every $i\in J$. Hence $V(\sum_{j\in J}I_{j})\subseteq \bigcap_{j\in J}V(I_{j})$. If $P\in \bigcap_{j\in J}V(I_{j})$, then $I_{i}\subseteq P$, for every $i\in J$. So $I_{i}\subseteq \sum_{j\in J}I_{j}\subseteq P$, and  $P\in V(\sum_{j\in J}I_{j})$. Therefore $V(\sum_{j\in J}I_{j})\subseteq \bigcap_{j\in J}V(I_{j})$.\\
		Now let $I$ and $J$ be hyperideals of $R$ and $P\in V(IJ)$. Since $P$ is prime and $IJ\subseteq P$, then $I\subseteq P$ or $J\subseteq P$. So $P\in V(I)\cup V(J)$. Let $P\in V(I)\cup V(J)$, then $I\subseteq P$ or $J\subseteq P$, so $IJ\subseteq P$ and $P\in V(IJ)$. Therefore $V(IJ)=V(I)\cup V(J)$.
	\end{proof}
\end{Theorem}

The resulting topology on $Spec(R)$ is called the Zariski topology of $R$, and the family $\{V_{m}(I)\}_{I\lhd R}$ forms a subspace for it.

\begin{remark}    \label{iii-6}
	Let $\{S_{j}\}_{j\in J}$ be a family of subsets of hyperring $R$. Since
	\begin{equation}
		I(\bigcup_{j\in J} S_{j})=\sum_{j\in J}I(S_{j}),
	\end{equation}
	then by Lemma \ref{iii-4},  for every $i,j\in J$ we have
	\begin{center}	
		$V(\bigcup_{j\in J} S_{j})=\bigcap_{j\in J}V(S_{j})=\bigcap_{j\in J}V(I(S_{j}))=V(\sum_{j\in J} I(S_{j}))$,
	\end{center}
	\begin{center}
		$V(S_{i}\cap S_{j})=V(S_{i})\cup V(S_{j})=V(I(S_{i}))\cup V(I(S_{j}))=V(I(S_{i})\cap I(S_{j}))$.
	\end{center}
	Also, if $I$ and $J$ are hyperideals of $R$ then by Lemma \ref{iii-4}, $V(I)\cup V(J)\subseteq V(I\cap J)$. Let $P\in V(I\cap J)$ but $P\notin V(I)\cup V(J)$, then there are $x\in I-P$ and $y\in J-P$ such that $xy\in I\cap J$. So $xy\in P$, and it is a contradiction. Therefore $V(IJ)=V(I\cap J)$.
\end{remark}

\begin{proposition}     \label{iii-7}
	If $I$ and $J$ are hyperideals of $R$, then:
		\item[(i)] {$V(I)=V(\sqrt{I})$};
		\item[(ii)] { $V(I)\subseteq V(J)\Leftrightarrow\sqrt{J}\subseteq \sqrt{I}$}.
	\begin{proof}
		It will be proven similar to classic rings.
	\end{proof}
\end{proposition}

\begin{Theorem}      \label{iii-8}
	If $f:R\rightarrow S$ is a good homomorphism of hyperrings, then $\bar{f}:Spec(S)\rightarrow Spec(R)$, defined by $\bar{f}(P)=f^{-1}(P)$, is continuous.
	\begin{proof}
		Let $V(I)$ be a closed set of $Spec(R)$. Then:
		\begin{align*}
			\bar{f}^{ ^{ ^{-1}}}(V(I))&=\{P\in Spec(S) ; \ \bar{f}(P)\in V(I)\}=\{P\in Spec(S) ; \ I\subseteq f^{-1}(P) \}\\
			&=\{P\in Spec(S) ; \ f(I)\subseteq P \}=V(f(I)).
		\end{align*}
	\end{proof}
\end{Theorem}

If $R$ is a hyperring and $I$ is a hyperideal of $R$, then $\pi:R\rightarrow R/I$ is projection map and $\bar{\pi}$ is a continuous map from $Spec(R/I)$ to $Spec(R)$.

\begin{Theorem}     \label{iii-9}
	Let $R$ be a hyperring and $I$ be a hyperideal of $R$. Then:
		\item[(i)] {$\bar{\pi}(Spec(R/I))=V(I)$};
		\item[(ii)] {$\bar{\pi}$ is injective};
		\item[(iii)] {$Spec(R/I)$ is homiomorphic to $V(I)$ by subspace topology}.
	\begin{proof}
		We have
		\begin{align*}
			\bar{\pi}(Spec(R/I))&=\{P\in Spec(R); \ P/I\in Spec(R/I) \} \\
			&=\{P\in Spec(R); \ I\subseteq P \}\\
			&=V(I).
		\end{align*}
		Since there is a bijection between the hyperideals of $R/I$ and the hyperideals of $R$ containing $I$, then $\bar{\pi}$ is a bijection between $Spec(R/I)$ and $V(I)$.\\
		By Theorem \ref{iii-8} $\bar{\pi}$ is continuous and also $\bar{\pi}^{^{-1}}:V(I)\rightarrow Spec(R/I)$ defined by $\bar{\pi}^{^{-1}}(P)=P/I$ is continuous, because every prime hyperideal of $R/I$ is of the form $P/I$, where $P\in V(I)$.
	\end{proof}
\end{Theorem}

\begin{corollary}   \label{iii-10}
	If $R$ is a hyperring, then $Spec(R)=Spec(R/nil(R))$.
	\begin{proof}
		By Theorem \ref{iii-9}, $Spec(R/nil(R))=V(nil(R))$ and since $nil(R)\subseteq P$, for every $P\in Spec(R)$ we have $V(nil(R))=Spec(R)$.
	\end{proof}
\end{corollary}

\begin{Theorem}     \label{iii-11}
	Let $R_{1}$ and $R_{2}$ be hyperrings, then:
	\begin{equation}	
		Spec(R_{1}\times R_{2})=Spec(R_{1}) \ \dot{\cup} \ Spec(R_{2}).
	\end{equation}
	\begin{proof}
		Let $P_{1}\times P_{2}$ be a prime hyperideal of $R_{1}\times R_{2}$. Then $R_{1}/P_{1}\times R_{2}/P_{2}$, is hyperdomain. So $P_{1}=R_{1}$ and $P_{2}\in Spec(R_{2})$ or $P_{2}=R_{2}$ and $P_{1}\in Spec(R_{1})$.
	\end{proof}
\end{Theorem}

\begin{definition}[\cite{bib13}]     \label{iii-12}
	The hyperideals $I$ and $J$ of hyperring $R$ said to be comaximal if $I+J=R$.
\end{definition}

\begin{Theorem}[Chinese Remainder Theorem]     \label{iii-13}
	If $I_{1},I_{2},...,I_{k}$ are hyperrings of $R$, then the map $R\rightarrow R/I_{1}\times R/I_{2}\times...\times R/I_{k}$ defined by $r\longmapsto (r+I_{1},r+I_{2},...,r+I_{k})$ is a good homomorphism with kernel $I_{1}\cap I_{2}\cap...\cap I_{k}$. If for each $i,j\in \{1,2,...,k \}$ with $i\neq j$ the hyperideals $I_{i}$ and $I_{j}$ are comaximal, then this map is surjective and $I_{1}\cap I_{2}\cap...\cap I_{k}=I_{1} I_{2}...I_{k}$, so:
	\begin{equation}
		R/(I_{1} I_{2}...I_{k})=R/(I_{1} \cap I_{2}\cap...\cap I_{k})=R/I_{1}\times R/I_{2}\times...\times R/I_{k}.
	\end{equation}
	\begin{proof}
		Let $I=I_{1}$ and $J=I_{2}$. Consider the map $\varphi:R\rightarrow R/I \times R/J$, defined by $\varphi(r)=(r+I,r+J)$. Then $\varphi$ is good homomorphism of hyperrings and $Ker\varphi=I\cap J$. Since $R$ is a  Krasner hyperring with a unit
		element and $I+J=R$, then there are elements $x\in I$ and $y\in J$ such that $1\in x+y$. Hence $x\in1-y\subseteq1+J$ and
		$y\in1-x\subseteq1+I$. Also $1-y+J=1+J$ and $1-x+I=1+I$. So $\varphi(x)=(I,1+J)=(0_{R/I},1_{R/J})$ and $\varphi(y)=(1+I,J)=(1_{R/I},0_{R/J})$. Let $(r_{1}+I,r_{2}+J)\in R/I \times R/J$, then
		\begin{align*}
			(r_{1}+I,r_{2}+J)&=(r_{1}+I,0)+(0,r_{2}+J) \\
			&=(r_{1}+I,r_{1}+J)(1,0)+(r_{2}+I,r_{2}+J)(0,1)  \\
			&=\varphi(r_{1})\varphi(y)+\varphi(r_{2}) \varphi(x)  \\
			&=\varphi(r_{1}y+r_{2}x).
		\end{align*}
		We know that $IJ\subseteq I\cap J$. Let $z\in I\cap J$, then $z=z.1\in z(x+y)=zx+zy\in IJ$. Hence $I\cap J\subseteq IJ$. The general case will follow by induction from the case of two hyperideals using $I=I_{1}$ and $J=I_{2}I_{3}...I_{K}$.
	\end{proof}
\end{Theorem}

\begin{corollary}      \label{iii-14}
	The hyperring $R$ is a product of hyperrings if and only if $R$ has nontrivial idempotents.
	\begin{proof}
		Let $x$ be a nontrivial idempotent of $R$. Hence $(x)$ and $(1-x)$ are comaximal hyperideals of $R$ and $(x)(1-x)=(0)$. Now by the Chinese Remainder Theorem: $R\cong R/(0) \cong R/(x)\times R/(1-x)$. Conversely suppose that $R=R_{1}\times R_{2}$, then $(1,0)^{2}=(1,0)$.
	\end{proof}
\end{corollary}

\begin{lemma}     \label{iii-15}
	If $R/nil(R)$ has nontrivial idempotent then $R$ has nontrivial idempotent.
	\begin{proof}
		Let $x+nil(R)\in R/nil(R)$ be nontrivial idempotent. So $x^{2}+nil(R)=x+nil(R)$ and $x^{2}-x\subseteq nil(R)$.
		Hence there is $n\in \mathbb{N}$ such that $(x^{2}-x)^{n}=x^{n}(x-1)^{n}=0$. Also $(x^{n})+((x-1)^{n})=R$, and $(x^{n}).((x-1)^{n})=(0)$.
		Now by Chinese Remainder Theorem $\varphi:R\rightarrow R/(x^{n})\times R/((x-1)^{n})$, defined by $\varphi(x)=\Big(r+(x^{n}),r+((x-1)^{n})\Big)$ is an isomorphism and $(0,1)\in R/(x^{n})\times R/((x-1)^{n})$ is nontrivial idempotent. So $\varphi^{-1}((0,1))\in R$ is nontrivial idempotent.
	\end{proof}
\end{lemma}

\begin{Theorem}     \label{iii-16}
	If $Spec(R)$ is a disconnected space and $Spec(R)=C \ \dot{\cup}\ D$, then $R\cong R_{1}\times R_{2}$, that $C=Spec(R_{1})$ and $D=Spec(R_{2})$.
	\begin{proof}
		Let $I$ and $J$ be hyperideals of $R$ such that $C=V(I)$ and $D=V(J)$. Since $V(R)=\emptyset$ and $V(0)=Spec(R)$, then
		$$\emptyset=C\cap D=V(I)\cap V(J)=V(I+J)=V(R);$$
		$$Spec(R)=C\cup D=V(I)\cup V(J)=V(IJ)=V(0).$$
		By Lemma \ref{iii-4} we  have $I+J=R$ and by Chinese Remainder Theorem $R/IJ\cong R/I\times R/J$. If $R$ has no nilpotent elements, then $nil(R)=(0)$, and by Proposition \ref{iii-7} $IJ=0$. So $R/(0)\cong R\cong R/I\times R/J$, and by Theorem \ref{iii-9} we have $V(I)=Spec(R/I)$ and $V(J)=Spec(R/J)$.
		Now since $Spec(R)=Spec(R/nil(R))$ so $Spec(R/nil(R))$ is disconnected and $R/nil(R)$ has no nilpotent elements.
		Hence $R/nil(R)\cong S\times T$, for some hyperrings $S$ and $T$, and by Corollary \ref{iii-14}, $R/nil(R)$ has nontrivial idempotents.
		Therefore by Lemma \ref{iii-15}, $R$ has nontrivial idempotent and by Corollary \ref{iii-14}, $R$ is product of hyperrings.
	\end{proof}
\end{Theorem}

\begin{remark}     \label{iii-17}
	Consider $W(f)=\big\{P\in Spec(R); f\notin P \big\}$, for every $f\in R$. Then $W(0)=\emptyset$, $W(1)=X=Spec(R)$ and $X-V(E)=\bigcup_{f\in E}W(f)$, for $E\subseteq R$. Therefore $\{W(f); f\in R \}$ is a basis for Zariski topology on $Spec(R)$, and it is clear that $W(f)\cap W(g)=W(fg)$. Also $W(f)=W(g)$ if and only if $\sqrt{(f)}=\sqrt{(g)}$, and $f\in R$ is nilpotent if and only if $W(f)=\emptyset$, and $f\in R$ is unit if and only if $W(f)=X$. Here too like the theory of classical rings, it can be proved that $W(f)$ is quasi-compact for every $f\in R$, and an open subset of $X$ is quasi-compact if and only if it is a finit union of sets $W(f)$ \cite{bib5}.
	On the other hand if we define Zariski topology based on closed subsets, then $\mathbb{B}=\{B(x); x\in R\}$ is a basis for Zariski topology on $Spec(R)$, where $B(x)=V(( x))=\{P\in Spec(R);  x\in P\}$.
\end{remark}

\begin{Theorem}     \label{iii-18}
	The hyperideal $nil(R)$ of $R$ is prime if and only if $Spec(R)$ is irreducible, (especially that if $R$ is hyperdomain, then $Spec(R)$ is irreducible).
	\begin{proof}
		We know that $\{W(f) \}_{f\in R}$, is a basis for Zariski topology on $Spec(R)$ and any two non-empty sets will intersect if and only if any tow non-empty basis elements intersect.
		So $Spec(R)$ is irreducible if and only if any tow non-empty basis elements intersect, that is for every $W(f)\neq \emptyset$ and $W(g)\neq \emptyset$ we have $W(fg)=W(f)\cap W(g)\neq \emptyset$. Therefore $Spec(R)$ is irreducible if and only if for every $f,g\in R$, if $f,g\notin nil(R)$ then  $fg\notin nil(R)$, which means that $nil(R)$ is the prime hyperideal of $R$.
	\end{proof}
\end{Theorem}

\begin{corollary}      \label{iii-19}
	Let $R$ be a hyperring and $I$ be a hyperideal of $R$. Then $V(I)$ is an irreducible subset of $Spec(R/I)$ if and only if $\sqrt{I}$ is prime.
	\begin{proof}
		By Theorem \ref{iii-18}, $Spec(R/I)$ is irreducible if and only if $\sqrt{I}$ is prime. Also $Spec(R/I)\cong V(I)$.
	\end{proof}
\end{corollary}

\begin{remark}      \label{iii-20}
	The subset $\{x\}=\{P \}$ of $Spec(R)$ is closed (closed point), if and only if $P\in mSpec(R)$. So $\overline{\{x\}}=V(P)$ and $y=Q\in \overline{\{x\}}$ if and only if $P\subseteq Q$. Also $Spec(R)$ is $T_{0}$-Space, because if $x=P$ and $y=Q$ are distinct points of $Spec(R)$, then $P\nsubseteq Q$ or $Q\nsubseteq P$. Without loss of generality assume that $Q\nsubseteq P$.
	Then $x \notin \overline{\{y\}}$ and so $Spec(R)-\overline{\{y\}}$ is an open set that contains $x$, and $y\notin Spec(R)-\overline{\{y\}}$.
\end{remark}

\begin{definition}[\cite{bib20}]      \label{iii-21}
	If $P\in Spec(R)$, the height of $P$ denoted by $h(P)$ and defined to be the supremum of lengths of chains $P_{0}\subsetneq P_{1}\subsetneq...\subsetneq P_{n}=P$, if this supremum exists, and $\infty$ otherwise. The dimension of $R$, denoted by $dim(R)$ and $dim(R)=sup\{h(P); P\in Spec(R)\}$. It is clear that if $R$ is a hyperfield then $dim(R)=0$.
\end{definition}

\begin{Theorem}      \label{iii-22}
	$Spec(R)$ is $T_{1}$-Space if and only if $dim(R)=0$ if and only if $Spec(R)$ is Hausdorff space.
	\begin{proof}
		If $dim(R)\neq0$, then there are $P_{1},P_{2}\in Spec(R)$ such that $P_{1}\subsetneq P_{2}$. So every neighborhood containing $P_{2}$ is contain $P_{1}$, and it is contradiction. Now suppose that $dim(R)=0$. So there are $P_{1},P_{2}\in Spec(R)$ such that $x\in P_{1}- P_{2}$ and $y\in P_{2}- P_{1}$. Therefore $(x)\subseteq P_{1}$,  $(x)\nsubseteq P_{2}$ and $(y)\subseteq P_{2}$,  $(y)\nsubseteq P_{1}$. So $P_{1}\in Spec(R)-V((y))$, $P_{1}\notin Spec(R)-V((x))$ and $P_{2}\in Spec(R)-V((x))$, $P_{2}\notin Spec(R)-V((y))$. Every Hausdorff space ($T_{2}$-Space), is $T_{1}$-Space. So assume that $dimR=0$, and $P,Q\in Spec(R)$ and $x\in P-Q$. We know that $P_{P}$ is maximal hyperideal of local hyperring $R_{P}$, and $P_{P}=nil(R_{P})$. Let $\frac{x}{1}\in P_{P}=nil(R_{P})$. Then there is $n\in\mathbb{N}$ such that $\frac{x^{n}}{1}=0_{R_{P}}$. Hence for some $s\in R-P$ we have $sx^{n}=0$.\\
		Since $x\notin Q$ and $s\notin P$ then $Q\in Spec(R)-V((x))=W((x))$ and $P\in Spec(R)-V((s))=W((s))$. Also $W((x))\cap W((s))=\{I\in Spec(R); s\notin I, x\notin I \}=\{I\in Spec(R); s\notin I, x^{n}\notin I  \}=\{I\in Spec(R); sx^{n}\notin I \}=Spec(R)-V((sx^{n}))=\emptyset$, since $V(( sx^{n}))=V((0))=Spec(R)$.
	\end{proof}
\end{Theorem}

\begin{corollary}          \label{iii-23}
	$Spec(P)$ is $T_{1}$-space if and only if $mSpec(R)=Spec(R)$.
\end{corollary}

\begin{proposition}           \label{iii-24}
	For any hyperring $R$, $Spec(R)$ is a compact space under the Zariski topology.
	\begin{proof}
		Let $\{I_{j}\}_{j\in J}$ be a family of hyperideals of $R$ and $\{W(I_{j})\}_{j\in J}$ be a family of open sets that $Spec(R)=\bigcup_{j\in J}W(I_{j})$ and $W(I_{j})=Spec(R)-V(I_{j})$, for any $j\in J$. So $Spec(R)=\bigcup_{j\in J}(Spec(R)-V(I_{j}))=Spec(R)-\bigcap_{j\in J}V(I_{j})=Spec(R)-V(\sum_{j\in J}I_{j})$. Hence $V(\sum_{j\in J}I_{j})=\emptyset$ and $\sum_{j\in J}I_{j}=R$. So $1\in\sum_{k\in K}a_{k}x_{k}$ where $a_{k}\in R$, $x_{k}\in I_{k}$ and $K$ is a finite subset of $J$. Therefore $\sum_{k\in K}I_{k}=R$ and $V(\sum_{k\in K}I_{k})=\emptyset$. So $Spec(R)=\bigcup_{k\in K}W(I_{k})$.
	\end{proof}
\end{proposition}

\begin{proposition}       \label{iii-25}
	Let $R$ be hyperring and $x\in R$. Then $W(x)$ is a compact subset of $Spec(R)$. Also an open subset of $Spec(R)$ is compact if and only if it is a finite union of sets $W(x_{i})$, where $x_{i}\in R$.
	\begin{proof}
		Consider $W(x)=\bigcup_{i\in I}W(x_{i})$. For every $P\in Spec(R)$, if $x\notin P$ then there is $i\in I$ such that $x_{i}\notin P$. Let $K=\sum_{i\in I}(x_{i})$. Since $Spec(R)-V(x)=Spec(R)-V(K)$, then for every prime hyperideal $P$, that $K\subseteq P$, we have $x\in P$. So $x\in \sqrt{K}$ and $x^{m}\in K$, for some $m\in \mathbb{N}$. Hence for some $x_{1},x_{2},...,x_{n}$ and $r_{i}\in R  (1\leq i\leq n)$, we have $x^{m}=\sum_{i=1}^{n}r_{i}x_{i}$. If $P\in Spec(R)$ and $\{x_{1},x_{2},...,x_{n}\}\subseteq P$, then $x^{m}\in P$ and so $x\in P$. Also if $x\notin P$ then $x_{i}\notin P$, for some $i\in\{1,2,...,n\}$. Therefore $W(x)=\bigcup_{i=1}^{n}W(x_{i})$.\\
		For the last part, let $X$ be a compact and open subset of $Spec(R)$ and consider $X=\bigcup_{i\in I}W(x_{i})$. Then there is a finite subcover $X=\bigcup_{i=1}^{n}W(x_{i})$ for some $x_{1},x_{2},...,x_{n}$. Also $X=\bigcup_{i=1}^{n}W(x_{i})$ is compact, since the finite union of compact sets is compact.
	\end{proof}
\end{proposition}

\begin{remark}[\cite{bib5}]          \label{iii-26}
	For the hyperring $R$, also like the rings, it can be shown that if $R$ is Arthinian, then every prime hyperideal is maximal, and the number of maximal hyperideals is finite. Also, $R$ is Arthinian if and only if $R$ is Noetherian and $dim(R)=0$.
\end{remark}

By using Proposition \ref{iii-7} it can be easily proved that if $R$ is a Noetherian hyperring, then $Spec(R)$ is Noetherian space.

\begin{proposition}      \label{iii-28}
	Let $R$ be a Noetherian hyperring. The following assertions are equivalent:
		\item[(i)] {$R$ is Arthinian};
		\item[(ii)] {$Spec(R)$ is discrete and finite};
		\item[(iii)] {$Spec(R)$ is discrete}.
	\begin{proof}
		$1\Rightarrow2$: If $R$ is Arthinian hyperring by Remark \ref{iii-26}, $mSpec(R)=Spec(R)$ and each point of $Spec(R)$ is closed and the number of maximal hyperideals of $R$ is finite. Hence $Spec(R)$ is discrete and finite.\\
		$3\Rightarrow1$: Since $Spec(R)$ is discrete, then each point is closed. By Lemma \ref{iii-4}, any closed point is maximal hyperideal. So $mSpec(R)=Spec(R)$ and $dim(R)=0$. Now by Remark \ref{iii-26}, $R$ is Arthinian.
		
	\end{proof}
\end{proposition}

\begin{proposition}      \label{iii-29}
	Let $R$ be hyprring. The irreducible components of $Spec(R)$, are closed sets $V(I)$, where $I$ is a minimal prime hyperideal of $R$.
	\begin{proof}
		If $X$ is a maximal irreducible subset of $Spec(R)$, then by Theorem \ref{ii-13}, $X$ is closed and so $X=V(I)$, for some hyperideal $I$ of $R$. By Corollary \ref{iii-19}, $\sqrt{I}$ is prime hyperideal and if $P\in Spec(R)$ such that $V(P)$ is irreducible and $P\subseteq I$, then $X=V(I)\subseteq V(P)$. So $V(I)= V(P)$ and therefore $\sqrt{P}=\sqrt{I}$.
	\end{proof}
\end{proposition}

\begin{Theorem}        \label{iii-30}
	If $R$ is a local hyprring, then $Spec(R)$ is connected.
	\begin{proof}
		Let $x\in R$ be a nontrivial idempotent and $M$ be a maximal hyperideal of $R$. Then $x\notin M=nil(R)$. Since $x$ is not unit, then $X$ is contained in some maximal hyperideal. So $x\in M$ and it is contradiction.
	\end{proof}
\end{Theorem}

\begin{remark}[\cite{bib9}]        \label{iii-32}
	Let $S$ be a multiplicative subset of hyperring $R$, and $I$ be a hyperideal of $R$. Then $S\cap I\neq \emptyset$ if and only if $S^{-1}I=S^{-1}R$. Also there is a one to one correspondence between the sets $\{P; P\in Spec(R), P\cap S=\emptyset\}$ and $Spec(S^{-1}R)$, under the mapping $P\mapsto S^{-1}P$.
\end{remark}

\begin{proposition}         \label{iii-33}
	Let $S$ be a multiplicative subset of hyperring $R$. Then $Spec(R)$ is irreducible and $nil(R)\cap S=\emptyset$ if and only if $Spec(S^{-1}R)$ is irreducible.
	\begin{proof}
		If $Spec(R)$ is irreducible, then by Theorem \ref{iii-18}, $nil(R)\in Spec(R)$ and since $nil(R)\cap S=\emptyset$ and $nil(S^{-1}R)=S^{-1}nil(R)$, then $nil(S^{-1}R)\in Spec(S^{-1}R)$ and so $Spec(S^{-1}R)$ is irreducible.\\
		If $Spec(S^{-1}R)$ is irreducible, then $nil(S^{-1}R)=S^{-1}nil(R)\in Spec(S^{-1}R)$. Hence $nil(R)\cap S=\emptyset$ and $nil(R)\in Spec(R)$.
	\end{proof}
\end{proposition}

\begin{proposition}        \label{iii-34}
	$Spec(S^{-1}R)$ is disconnected if and only if $Spec(R)$ is disconnected.
	\begin{proof}
		Since $S^{-1}(R_{1}\times R_{2})\cong S^{-1}R_{1}\times S^{-1}R_{2}$, by Theorem \ref{iii-11} and Theorem \ref{iii-16} the result is obtained.
	\end{proof}
\end{proposition}

\begin{example}       \label{iii-35}
	Let $(R,+,.)$ be a commutative ring with a unit element and $G$ be a subgroup of monoid $(R-\{0\},.)$ and $I$ be an ideal of $R$.\\
	(1) Consider $\bar{R}=R/G=\{rG; r\in R\}$ and $rG\oplus sG=\{ tG; \ t\in rG+sG\}$, $rG\odot sG=rsG$. Then $(\bar{R},\oplus,\odot)$ is hyperring and $0_{\bar{R}}=\{0\}$, $1_{\bar{R}}=G$. Moreover if $R$ is a field, then $(\bar{R},\oplus,\odot)$ is a hyperfield,(For more details see \cite{bib15}).\\
	The hyperideals of $R$ are of the form $I/G$ such that $I$ is ideal of $R$. Let $P\in Spec(R)$ and $rG\odot sG=rsG\in P/G$, then $rs\in P$. So $r\in P$ or $s\in P$ and therefore $rG\in P/G$ or $sG\in P/G$. and $P/G\in Spec(\bar{R})$. Also for any $A\in Spec(\bar{R})$, there is a prime ideal $P\in Spec(R)$ such that $A=P/G$. So $Spec(\bar{R})=\{P/G; P\in Spec(R)\}$.\\
	(2) Consider $\bar{R}=R/I=\{r+I; r\in R\}$ and define $r+I\oplus s+I=(r+s)+I$ and $r+I\odot s+I=\{t+I; t\in rs+I\}$. Then $(\bar{R},\oplus,\odot)$ is multiplicative hyperring and $0_{\bar{R}}=I$, $1_{\bar{R}}=1+I$. It is clear that $Spec(\bar{R})=\{P/I; P\in Spec(R)\}\cong V(I)$.\\
	(3) If $\bar{R}=\{r+I; r\in R\}$, $r+I\oplus s+I=\{t+I; t\in r+s+I\}$ and $r+I\odot s+I=\{t+I; t\in rs+I\}$, then $(\bar{R},\oplus,\odot)$ is general hyperring and $Spec(\bar{R})=\{P/I; P\in Spec(R)\}\cong V(I)$.
\end{example}

\begin{corollary}      \label{iii-36}
	$Spec(\bar{R})$ is irreducible (disconnected) if and only if $Spec(R)$ is irreducible (disconnected).
	\begin{proof}
		We know that $P\in Spec(R)$ if and only if $P/G\in Spec(\bar{R})$. Since $1_{R}\in G$, then $nil(\bar{R})=\{rG;\ r^{m}G=0,\ for \ some \ m\in \mathbb{N}\}=\{rG;\ r^{m}=0,\ for \ some \ m\in \mathbb{N}\}=\{rG;\ r\in nil(R)\}=nil(R)/G$. So $nil(R)\in Spec(R)$ if and only if $nil(\bar{R})=nil(R)/G\in Spec(\bar{R})$. By Theorem \ref{iii-18}, $Spec(R)$ is irreducible if and only if $Spec(\bar{R})$ is irreducible.\\
		Since $(R_{1}\times R_{2})/G \cong R_{1}/G \times R_{2}/G$, it is clear that $Spec(R)$ is disconnected if and only if $Spec(\bar{R})$ is disconnected.
	\end{proof}
\end{corollary}

\begin{example}        \label{iii-37}
	Let $R=\mathbb{Z}$ and $G=\{\pm1\}$. Then $\bar{R}=\{ \{\pm a\}; a\in \mathbb{Z}\}$ and $\{\pm a\}\oplus \{\pm b\}=\{\{\pm (a+b)\}, \{\pm (a-b)\}\}$, $\{\pm a\}\odot \{\pm b\}=\{\pm ab\}$.
	For every prime $p\in \mathbb{Z}$; $(p)/G=\{\{\pm kp\}; k\in \mathbb{N}\cup \{0\} \}$. Also $\mathbb{Z}$ is integral domain therefore $\bar{R}$ is hyperdomain and $nil(\bar{R})$ is prime hyperideal of $\bar{R}$. So $Spec(\bar{R})$ is irreducible and connected.
\end{example}

\begin{example}        \label{iii-38}
	In Example \ref{iii-37}, if $F$ is a field and $R=F[x]$ and $G=F-\{0\}$, then $Spec(\bar{R})=\big\{(f(x))/G \ ; \ f(x)\in F[x] \ is\ irreducible \big\}$.
\end{example}

\section{From Zariski topology of hyperrings to Zariski topology of rings} \label{sec4}

\begin{definition}[\cite{bib10}]   \label{iv-1}
	Let $\mathcal{U}$ be the set of all finite sums of products of elements of general hyperring $R$. Then $(a,b)\in \gamma^{*}$ if and only if there exists $(a=z_{1},z_{2},...,b=z_{n+1})\in R^{n+1}$ and $U_{1},U_{2},...,U_{n}\in \mathcal{U}$ such that $\{z_{i},z_{i+1}\}\subseteq U_{i}$, for any $i\in \{1,2,...,n\}$.
\end{definition}

The relation $\gamma^{*}$ is the smallest equivalence relation on general hyperring $R$ such that the quotient $R/\gamma^{*}$ is a ring. $R/\gamma^{*}$ is called the fundamental ring, \cite{bib11}.

\begin{proposition}
	If $\rho$ is a strongly regular relation on Krasner hyperring $R$, then $\rho(0)=\{x\in R; \rho(x)=0_{R/\rho}\}$ is a normal hyperideal of $R$.
	\begin{proof}
		For every $x,y\in \rho(0)$;\ $\rho(x-y)=\rho(x)-\rho(y)=0_{R/\rho}$. So $x-y\subseteq\rho(0)$. Let $r\in R$, then $\rho(rx)=\rho(r)\rho(x)=\rho(r)0_{R/\rho}=0_{R/\rho}$. So $rx\in\rho(0)$, and similarly $xr\in\rho(0)$. Also, since $\gamma^{\ast}\subseteq\rho$ and $r-r\subseteq\gamma^{\ast}(0)$, then $r-r\subseteq\rho(0)$ and $r+\rho(0)-r\subseteq\rho(0)$. Thus $\rho(0)$ is normal.
	\end{proof}
\end{proposition}

Clearly every hyperideal $I$ of $R$, containing $\gamma^{\ast}(0)$ is normal, and $R/I$ is Krasner hyperring. Let us denote the set of all strongly regular relations on $R$, by $\mathcal{SR}(R)$ and the set of all hyperideals containing $\gamma^{\ast}(0)$, by $\mathcal{I}(\gamma^{\ast}(0))$.

\begin{lemma}   \label{22A} 
	If $\rho\in\mathcal{SR}(R)$, then $\rho(x)=x+\rho(0)$, for every $x\in R$.
	\begin{proof}
		If $z\in \rho(x)$, then $-x+z\subseteq\rho(0)$ and $z\in x+(-x+z)\subseteq x+\rho(0)$. So $z\in x+\rho(0)$ and $\rho(x)\subseteq x+\rho(0)$. Conversely, if $z\in x+\rho(0)$, then $z\in x+y$, for some $y\in\rho(0)$. So $\rho(z)=\rho(x)$ and $z\in \rho(x)$. Therefore $x+\rho(0)\subseteq\rho(x)$.	
	\end{proof}
\end{lemma}

\begin{lemma}   \label{22B} 
	If $I\in\mathcal{I}(\gamma^{\ast}(0))$, then the congruence relation modulo $I$, is strongly regular.
	\begin{proof}
		Let $x,y,z \in R$ and $x+I=y+I$. Then $(z+x)+I=(z+y)+I$	and since $S_{\beta}\subseteq I$, then for every $r\in z+x$ and $s\in z+y$; $(z+x)+I=r+I$ and $(z+y)+I=s+I$. Hence the congruence relation modulo $I$ is strongly regular.		
	\end{proof}
\end{lemma}

Let $\mathcal{SR}(R)$ be the set of all strongly regular relations on canonical hypergroup $(R,+)$. 
If $\rho,\sigma\in \mathcal{SR}(R)$, then $\rho\vee\sigma\in \mathcal{SR}(R)$, and the lattices  $N(S_{\beta})=\{H;\ S_{\beta}\subseteq H\lhd R\}$ and $\mathcal{SR}(R)$ are isomorph, \cite{bib2}.

\begin{Theorem}
	For every hyperring $R$, $\mathcal{SR}(R)\cong\mathcal{I}(\gamma^{\ast}(0))$ is an isomorphism of complete lattices.
	\begin{proof}
		Since $R$ is a Krasner hyperring, then for every $\rho,\sigma\in\mathcal{SR}(R)$, $A\subseteq\mathcal{SR}(R)$, $I,J\in \mathcal{I}(\gamma^{\ast}(0))$ and $B\subseteq\mathcal{I}(\gamma^{\ast}(0))$, we have $\rho\vee\sigma,\rho\cap\sigma\in\mathcal{SR}(R)$ and $I\vee J=I+J, I\cap J \in \mathcal{I}(\gamma^{\ast}(0))$. Also
		\begin{center}
			$\underset{\rho\in A}{\bigvee}\rho\in\mathcal{SR}(R) ,\  \underset{\rho\in \mathcal{SR}(R)}{\bigcap}\rho=\gamma^{\ast},\  \underset{I\in B}{\bigvee}I=\underset{I\in B}{\sum}I\in \mathcal{I}(\gamma^{\ast}(0)), \ \underset{I\in \mathcal{I}(\gamma^{\ast}(0))}{\bigcap}I=\gamma^{\ast}(0)$.
		\end{center}
		Let $f:\mathcal{SR}(R)\rightarrow\mathcal{I}(\gamma^{\ast}(0))$ by $\rho\mapsto\rho(0)$ and $g:\mathcal{I}(\gamma^{\ast}(0))\rightarrow\mathcal{SR}(R)$ by $I\mapsto\rho_{I}:=\{(x,y)\in R^{2};\ x+I=y+I\}$. By Lemmas \ref{22A} and \ref{22B}, $f$ and $g$ are well defined, and
		\begin{align*}
			f\circ g(I)&=f(\rho_{I})\\
			&=\{r\in R; \ (x+I)+(r+I)=x+I, \forall x\in R\}\\
			&=\{r\in R; \ (x+r)+I=x+I, \forall x\in R\}=I,
		\end{align*}
		because, for every $x\in R$, $0\in-x+x\subseteq I$. Also 
		\begin{align*}
			g\circ f(\rho)&=g(\rho(0))\\
			&=\{(x,y)\in R^{2}; \ x+\rho(0)=y+\rho(0)\}\\
			&=\{(x,y)\in R^{2}; \ \rho(x)=\rho(y)\}=\rho.
		\end{align*}
		Now let $\rho,\sigma\in\mathcal{SR}(R)$ and $I,J\in\mathcal{I}(\gamma^{\ast}(0))$ such that $\rho\subseteq\sigma$ and $I\subseteq J$. Then $f(\rho)=\rho(0)\subseteq\sigma(0)=f(\sigma)$ and if $(x,y)\in R^{2}$ such that $x+I=y+I$, then $x+I+J=y+I+J$. Therefore $x+J=y+J$ and $f^{-1}(I)\subseteq f^{-1}(J)$.
	\end{proof}
\end{Theorem}

\begin{corollary}
	If $\rho\in\mathcal{SR}(R)$, then $R/\rho\cong R/\rho(0)$ is a ring isomorphism.
	\begin{proof}
		Since for every $x,y\in R$ and $z\in x+y$;\ $(x+\rho(0))+(y+\rho(0))=(x+y)+\rho(0)=z+\rho(0)$, then $R/\rho(0)$ is a ring. Consider the bijection map $\phi:R/\rho\rightarrow R/\rho(0)$ by $\rho(x)\mapsto x+\rho(0)$, then:
		$\phi(\rho(x)+\rho(y))=\phi(\rho(x+y))=(x+y)+\rho(0)=(x+\rho(0))+(y+\rho(0))=\phi(\rho(x))+\phi(\rho(y))$ and
		$\phi(\rho(x)\rho(y))=\phi(\rho(xy))=(xy)+\rho(0)=(x+\rho(0))(y+\rho(0))=\phi(\rho(x))\phi(\rho(y))$.	
	\end{proof}	
\end{corollary}

Therefore if $\rho,\sigma\in\mathcal{SR}(R)$ and $I,J\in \mathcal{I}(\gamma^{\ast}(0))$, then $\rho(0)+\sigma(0)=(\rho\vee\sigma)(0)$, $\rho(0)\cap\sigma(0)=(\rho\cap\sigma)(0)$ and $\rho_{I}\vee\rho_{J}=\rho_{I+J}$, $\rho_{I}\cap\rho_{J}=\rho_{I\cap J}$.\\
The hyperideal $I$ of a Krasner hyperring $R$ is a normal hyperideal if and only if $x+I-x\subseteq I$, for all $x\in R$, \cite{bib12}. Also it is proved that $R/I$ is a ring if and only if $I$ is a normal hyperideal of $R$, \cite{bib12}.

\begin{corollary}
	Let $I$ be a hyperideal of Krasner hyperring $R$. Then, $I$ is normal if and only if $\gamma^{\ast}(0)\subseteq I$.	
	\begin{proof}
		Since $0\in I$, then the hyperideal $I$ is normal if and only if $x-x\subseteq I$, for all $x\in R$. If $y\in \gamma^{\ast}(0)$, then there are $n\in \mathbb{N}$ and $(x_{1},...,x_{n})\in R^{n}$ such that $0,y\in \sum_{i=1}^{n}x_{i}$. Therefore $y\in \sum_{i=1}^{n}x_{i}-\sum_{i=1}^{n}x_{i}\subseteq I$. Let $\gamma^{\ast}(0)\subseteq I$. Since $x-x\subseteq \gamma^{\ast}(0)$, for all $x\in R$, then $I$ is normal.
	\end{proof}	
\end{corollary}

So we can say that if $R$ is a Krasner hyperring, then $R/I$ is ring if and only if $I\in \mathcal{I}(\gamma^{\ast}(0))$.

\begin{definition}
	A strongly regular relation $\rho$ on the commutative and of unity Krasner hyperring $R$ is called prime, if for every $x$ and $y$ in $R$, we have:
	\begin{equation}
		(xy,0)\in \rho \ \Rightarrow \ (x,0)\in\rho \ \ or \ \ (y,0)\in\rho.
	\end{equation}
	Also $\rho$ is primitive, if for every $x$ and $y$ in $R$, we have:
	\begin{equation}
		(xy,0)\in \rho \ , \ (x,0)\notin\rho \ \Rightarrow \ \exists n\in \mathbb{N} \ \ s.t.\ \ (y^{n},0)\in\rho.
	\end{equation}
\end{definition}

Since, $(x,0)\in\rho$ if and only if $x\in \rho(0)$, then $\rho$ is prime if and only if $xy\in\rho(0)$ results in $x\in\rho(0)$ or $y\in\rho(0)$, if and only if $\rho(x)\rho(y)=\rho(0)$ results in $\rho(x)=\rho(0)$ or $\rho(y)=\rho(0)$ if and only if $R/\rho$ is an integral domain. Therefore $\rho$ is prime if and only if $\rho(0)$ is a prime hyperideal of $R$.\\
Additionally, since the lattices $\mathcal{SR}(R)$ and $\mathcal{I}(\gamma^{\ast}(0))$ are isomorph, then $\rho$ is maximal if and only if $\rho(0)$ is maximal hyperideal of $R$.

\begin{example}
	Consider the Krasner hyperring $R$ as follows:
	\begin{center}
		\begin{tabular}{c|cccccccc}		
			$+$& $0$  & $1$      &$a$  &$b$  &$c$  &$d$ &$e$ & $f$\\
			\hline
			$0$& $0$  & $1$        &$a$       &$b$      &$c$      &$d$    &$e$     & $f$\\
			$1$& $1$  & $a$,$d$    &$b$,$e$   &$c$      &$1$,$f$  &$e$    &$0$,$c$ & $a$\\
			$a$& $a$  & $b$,$e$    &$0$,$c$   &$1$      &$a$,$d$  &$c$    &$1$,$f$ & $e$\\
			$b$& $b$  & $c$        &$1$       &$d$      &$e$      &$f$    &$a$     & $0$\\
			$c$& $c$  & $1$,$f$    &$a$,$d$   &$e$      &$0$,$c$  &$a$    &$b$,$e$ & $1$\\
			$d$& $d$  & $e$        &$c$       &$f$      &$a$      &$0$    &$1$     & $b$\\
			$e$& $e$  & $0$,$c$    &$1$,$f$   &$a$      &$b$,$e$  &$1$    &$a$,$d$ & $c$\\
			$f$& $f$  & $a$        &$e$       &$0$      &$1$      &$b$    &$c$     & $d$
		\end{tabular}
	\end{center}	
	\begin{center}
		\begin{tabular}{c|cccccccc}		
			$.$ & $0$  & $1$  &$a$   &$b$   &$c$     &$d$  &$e$     &$f$\\
			\hline
			$0$& $0$  & $0$  &$0$   &$0$   &$0$     &$0$  &$0$     &$0$\\
			$1$& $0$  & $1$  &$a$   &$b$   &$c$     &$d$  &$e$     &$f$\\
			$a$& $0$  & $a$  &$c$   &$d$   &$c$     &$0$  &$a$     &$d$\\
			$b$& $0$  & $b$  &$d$   &$f$   &$0$     &$d$  &$f$     &$b$\\
			$c$& $0$  & $c$  &$c$   &$0$   &$c$     &$0$  &$c$     &$0$\\
			$d$& $0$  & $d$  &$0$   &$d$   &$0$     &$0$  &$d$     &$d$\\
			$e$& $0$  & $e$  &$a$   &$f$   &$c$     &$d$  &$1$     &$b$\\
			$f$& $0$  & $f$  &$d$   &$b$   &$0$     &$d$  &$b$     &$f$			
		\end{tabular}
	\end{center}	
	Then $\gamma^{\ast}=\{(a,d),(d,a),(b,e),(e,b),(1,f),(f,1),(0,c),(c,0)\}\cup\Delta_{R}$ is a primitive strongly regular relation on $R$ that is not prime because $(ad,0)\in\gamma^{\ast}$ but $(a,0),(d,0)\notin\gamma^{\ast}$. Equivalently, $\gamma^{\ast}(0)=\{0,c\}$ is a primitive hyperideal that is not prime. Also, 
	$R/\gamma^{\ast}=\{\{0,c\},\{1,f\},\{a,d\},\{b,e\}\}\cong \mathbb{Z}_{4}$, that is not integral domain, and
	$\mathcal{I}(\gamma^{\ast}(0))=\{\{0,c\},\{0,c,a,d\}\}$, where $\{0,c,a,d\}=V(\gamma^{\ast}(0))=\sqrt{\gamma^{\ast}(0)}$ is a maximal hyperideal of $R$ and
	\begin{center}
	$R/\sqrt{\gamma^{\ast}(0)}=\{\{0,c,a,d\},\{1,f,b,e\}\}\cong \mathbb{Z}_{2}$
	\end{center}
	 is a field.
	Let $I=\{0,d\}$ and $M=\{0,b,d,f\}$. Because $a,f\notin I$ while $af=d\in I$ and also because  $\gamma^{\ast}(0)\nsubseteq I$, then $I$ is neither prime nor normal. Also, $M$ is a maximal hyperideal where is not normal. So, $Spec(R)=\{M, \sqrt{\gamma^{\ast}(0)}\}$.

\end{example}

\begin{Theorem}       \label{iv-3}
	Let $R$ be a hyperring, then $Spec(R/\gamma^{*})=\{P/\gamma^{*}; P\in Spec(R)\}$.
	\begin{proof}
		Let $I$ be a prime ideal of $R/\gamma^{*}$ and $P=\{s\in R; \gamma^{*}(r)\in I\}$. So $I=P/\gamma^{*}$ and if $s,t\in P$  and $r,k\in R$, then $\gamma^{*}(s)-\gamma^{*}(t)=\gamma^{*}(s-t)\in I$ and $\gamma^{*}(r)\gamma^{*}(s)=\gamma^{*}(rs)\in I$. Hence $s-t\subseteq P$ and $rs\in P$. If $\gamma^{*}(r)\gamma^{*}(k)=\gamma^{*}(rk)\in I$, then $\gamma^{*}(r)\in I$ or $\gamma^{*}(k)\in I$. So if $rk\in P$, then $r\in P$ or $k\in P$.
	\end{proof}
\end{Theorem}

\begin{lemma}    \label{iv-4}
	Let $S\subset Spec(R)$, then $\bar{S}=V(\bigcap_{P\in S}P)$.
	\begin{proof}
		Let $I$ be the radical hyperideal such that $V(I)=\bar{S}$, then $I\subset \bigcap_{P\in S}P$. Since $V(\bigcap_{P\in S}P)$ is closed, then $S\subset V(\bigcap_{P\in S}P)$ and so $\bar{S}\subset V(\bigcap_{P\in S}P)$. Hence $V(I)\subset V(\bigcap_{P\in S}P)$ and by Proposition \ref{iii-7}, $\bigcap_{P\in S}P=I$. 
		
	\end{proof}
\end{lemma}

\begin{Theorem}    \label{iv-5}   
	Let $f:R\rightarrow S$ be a good homomorphism of hyperrings. Then:
	\item[(i)] {If $I$ is a hyperideal of $R$, then $\bar{f}^{ ^{ ^{-1}}}(V(I))=V(I^{e})$};
		\item[(ii)] {For every $x\in R$, $\bar{f}^{^{^{-1}}}(W(x))=W(f(x))$};
		\item[(iii)] {If $J$ is a hyperideal of $S$, then $\overline{{\bar{f}(V(J))}}=V(J^{c})$};
		
		\item[(iv)] {If $f$ is surjective, then $\bar{f}$ is a homiomorphism of $Spec(S)$ on to the closed subset $V(Ker(f))$ of $Spec(R)$};
		\item[(v)] {If $f$ is injective, then $\bar{f}(Spec(S))$ is dense in $Spec(R)$. In fact, the image $\bar{f}(Spec(S))$ is dense in $Spec(R)$ if and only if $Ker(f)\subset \sqrt{0}$};
		\item[(vi)] {Let $g:S\rightarrow T$ be another homomorphism of hyperrings, then $\overline{g\circ f}=\bar{f}\circ\bar{g}$}.
	\begin{proof}
		(1)	By Theorem \ref{iii-8} we have $\bar{f}^{ ^{ ^{-1}}}(V(I))=V(f(I))$ and since $(f(I))=I^{e}$, by Lemma \ref{iii-4}, $V(f(I))=V(I^{e})$.\\
		(2) For every $x\in R$ we have:
		\begin{center}
			$\bar{f}^{ ^{ ^{-1}}}(W(x))=\{P\subset S ; x\notin f^{-1}(P)\}=\{P\subset S ; f(x)\notin P\}=W(f(x)).$
		\end{center}	
		(3)	By Lemma \ref{iv-4}, $\overline{\bar{f}(V(J))}=V(\bigcap_{P\in \bar{f}(V(J))}P)$ but
		\begin{center}
			$\bigcap_{P\in \bar{f}(V(J))}P=\bigcap_{J\subset Q}f^{-1}(Q)=
			f^{-1}(\bigcap_{J\subset Q}Q)=f^{-1}(\sqrt{J})=\sqrt{f^{-1}(J)}.$
		\end{center}
		Now by Proposition \ref{iii-7} we have $V(\bigcap_{P\in \bar{f}(V(J))}P)=V(f^{-1}(J))$.\\
		(4) Since $f$ is surjective, consider $S=R/Ker(f)$. We know that there is an inclusion perserving one to one correspondence between prime hyperideals of $R/Ker(f)$ and prime hyperideals of $R$ containing $Ker(f)$. So $\bar{f}$ is a continuous bijection onto the closed subset $V(Ker(f))$, and $\bar{f}^{ ^{ ^{-1}}}$ is continuous because:
		\begin{align*}
			\bar{f}(I/Ker(f))=&\{P\in Spec(R); I/Ker(f)\subset P/Ker(f)\in Spec(R/Ker(f))\}\\
			=&V(I).
		\end{align*}
		(5) Since $\bar{f}(Spec(S))=\bar{f}(V(0))$, then by (3) we have $\overline{\bar{f}(Spec(S))}=V(Ker(f))$. So $\bar{f}(Spec(S))$ is dense in $Spec(R)$ if and only if $V(Ker(f))=Spec(R)$ if and only if $Ker(f)\subset \sqrt{0}$.\\
		(6) It is clear that $\overline{g\circ f}(P)=(g\circ f)^{-1}(P)=(f^{-1}\circ g^{-1})(P)=\bar{f}(\bar{g}(P))$.
	\end{proof}
\end{Theorem}

Let $R$ be a hyperring. By $ZTop.HRg$, we mean the category of Zariski topology of (krasner) hyperrings, which its objects are $X=Spec(R)$ and for $X=Spec(R)$ and $Y=Spec(S)$, $Hom(X,Y)$ is the set of continuous maps induced by hyperring homomorphisms, with usual combinations of functions. Also, $ZTop.Rg$ denotes the category of Zariski topology of all rings.

\begin{Theorem} \label{iv-6}
	Consider the mapping $F:ZTop.HRg\rightarrow ZTop.Rg$ such that
	$Spec(R)\mapsto  Spec(R/\gamma^{*})$ and
	\begin{center}
	$(f:Spec(R)\rightarrow Spec(S))\mapsto  (f^{*}:Spec(R/\gamma^{*})\rightarrow Spec(S/\gamma^{*}))$.
	\end{center}
	 Then:
		\item[(i)] $F$ is functor;	
		\item[(ii)] The following diagram is commutative ($\pi_{R}$ and $\pi_{S}$ are canonical projections):
	\begin{center}
		\begin{tikzcd}	
			Spec(R) \ar[r,"f"]    \ar[d,"\pi_{R}"']        & Spec(S) \ar[d,"\pi_{S}"] \\
			Spec(R/\gamma^{*})  \ar[r,"f^{*}"']  &  Spec(S/\gamma^{*})
		\end{tikzcd}
	\end{center}
	\begin{proof}
		(1) Closed subsets of $Spec(S/\gamma^{*})$ are of the form $V(J/\gamma^{*})$ where $J$ is
		hyperideal of $S$, and $V(J/\gamma^{*})=\{P/\gamma^{*}\in Spec(S/\gamma^{*}); p\in V(J)\}$.
		So $g^{*^{-1}}(V(J/\gamma^{*}))=\{Q/\gamma^{*}\in Spec(R/\gamma^{*}); Q\in g^{-1}(V(J))\}$.
		Since $g$ is continuous, then $g^{-1}(V(J))$ is closed subset of $Spec(R)$ and hence there is hyperideal $I$ of $R$ such that $g^{-1}(V(J))=V(I)$. Thus $g^{*^{-1}}(V(J/\gamma^{*}))=V(I/\gamma^{*})$. It is clear that $F(1_{Spec(R)})=1_{Spec(R/\gamma^{*})}$ and $(h\circ g)^{*}=h^{*}\circ g^{*}$.\\
		(2) let $P\in Spec(R)$, then $f^{*}\pi_{R}(P)=f^{*}(P/\gamma^{*})=f(P)/\gamma^{*}=\pi_{S}f(P)$.
	\end{proof}
\end{Theorem}

\begin{Theorem}\label{iv-7}
	The mapping $Spec(-):H.Rg\rightarrow ZTop.HRg$ that $R\mapsto Spec(R)$ and $f\mapsto\bar{f}$ is a contravariant functor.
	\begin{proof}
		By Theorem \ref{iii-8}, $\bar{f}$ is continuous and by Theorem \ref{iv-5}, $\overline{g\circ f}=\bar{f}\circ\bar{g}$.
		
	\end{proof}
\end{Theorem}

\begin{Theorem}\label{iv-8}
	The following diagram is commutative:
	\begin{center}
		\begin{tikzcd}	
			H.Rg \ar[r,"Spec(-)"]    \ar[d,"\gamma^{*}"']        & ZTop.HRg \ar[d,"F"] \\
			Rg  \ar[r,"Spec(-)"']  &  ZTop.Rg
		\end{tikzcd}
	\end{center}	
	\begin{proof}
		It is clear that $F(Spec(R))=Spec(\gamma^{*}(R))$. If $f:R\rightarrow S$ is a good homomorphism of hyperrings then we have $\bar{f}:Spec(R)\rightarrow Spec(S)$ by $P\mapsto f^{-1}(P)$ and $f^{*}:R/\gamma^{*}\rightarrow S/\gamma^{*}$ by $\gamma^{*}(r)\mapsto \gamma^{*}(f(r))$. So $\bar{f}^{ ^{ ^{*}}}=\bar{f^{ ^{ *}}}$, because $\bar{f}^{ ^{ ^{*}}}(P/\gamma^{*})=f^{-1}(P)/\gamma^{*}=f^{*^{-1}}(P/\gamma^{*})=\bar{f^{ ^{ *}}}(P/\gamma^{*})$.	
		
	\end{proof}
\end{Theorem}

\section{Conclusions}

If we denote the set of all prime strongly regular relations on the hyperring $R$ by $V_{R}$, then $\{V(\rho)\}_{\rho\in \mathcal{SR}(R)}$ where $V(\rho)=\{\tau\in V_{R};\ \rho\subseteq\tau\}$ and $V_{R}=V(\gamma^{\ast})$ is a topological space, and for every $\rho,\sigma\in\mathcal{SR}(R)$:
\begin{center}
	$V(\rho)\cap V(\sigma)=V(\rho\vee\sigma) \ ,\ V(\rho)\cup V(\sigma)=V(\rho\cap\sigma)$.
\end{center}
Also, the topological space $V_{R}$ is homeomorphic to the topological subspace $V(\gamma^{\ast}(0))$ of $Spec(R)$.\\
Let's assume $H.Rg$ is the category of commutative Krasner hyperrings with identity and $\theta=\{\theta_{R}\}_{R\in H.Rg}$ is a sequence of strongly regular relations such that for every $R\in H.Rg$; $\theta_{R}\in V_{R}$. In this case, $\theta$ is a functor from the category $H.Rg$ to the category of integral domains that also maps Zariski topologies of hyperrings to Zariski topologies of rings. One application of this functor is in the study of sheaves of Krasner hyperrings.\\
In this paper, we studied the Zariski topology of Krasner hyperrings and investigated how fundamental relations act on them. The most important application of the results obtained here, one is that it will help our future work on the Zariski topology of Krasner hypermodules and other hyperrings, And the other is that they are used in the study of sheaves of hypermodules and hyperrings.


\begin{thebibliography}{99}	
	
		
	
	\bibitem{bib1}   B. Afshar and R. Ameri, Strongly Regular Relations Derived from Fundamental Relation, Journal of Algebraic Hyperstructures and Logical Algebras 4 (2023) 123-130, https://doi.org/10.61838/kman.jahla.4.2.8.
	
	
	\bibitem{bib2}R. Ameri and B. Afshar, Strongly regular relations on regular hypergroups, Journal of Mahani Mathematical Research 14 (2025) 73--83,  https://doi.org/10.22103/jmmr.2024.23228.1612.
	
	
	\bibitem{bib3} R. Ameri and A. Kordi, Clean multiplicative hyperrings, Italian Journal of Pure and Applied Mathematics 35 (2015) 625-636.
	
	
	\bibitem{bib4} S. M. Anvariyeh, S. Mirvakili and B. Davvaz, $\theta^{\ast}$-relation on hypermodules
	and fundamental modules over commutative fundamental rings, Communications
	in Algebra® 36 (2008) 622–631, http://dx.doi.org/10.1080/00927870701724078.
	
	
	
	\bibitem{bib5}  M. F. Atiyah and I. G. Macdonald, Introduction to Commutative
	Algebra, Addison-Wesley (1969), http://dx.doi.org/10.1017/S0008439500031039.
	
	
	
	\bibitem{bib6} A. Connes and C. Consani, From monoids to hyperstructures: in
	search of an absolute arithmetic, Casimir Force, Casimir Operators and
	the Riemann Hypothesis, de Gruyter (2010) 147–198, http://dx.doi.org/10.1515/9783110226133.147.
	
	
	\bibitem{bib7}  A. Connes and C. Consani, The hyperring of adele classes, Journal
	of Number Theory 131 (2011) 159–194, http://dx.doi.org/10.1016/j.jnt.2010.09.001. 
	
	
	
	\bibitem{bib8} A. Connes and C. Consani, Universal thickening of the field of real numbers, In Advances in the Theory of Numbers, Springer (2015) 11–74, http://dx.doi.org/10.1007/978-1-4939-3201-6-2.
	
	
	
	
	\bibitem{bib9} P. Corsini, Prolegomena of hypergroup theory, Aviani editore (1993).
	
	
	
	
	\bibitem{bib10} P. Corsini and V. Leoreanu, Applications of hyperstructure
	theory, Springer Science \& Business Mediavolume 5 (2003), http://dx.doi.org/10.1007/978-1-4757-3714-1. 
	
	
	\bibitem{bib11} B. Davvaz and V. Leoreanu-Fotea, Hyperring theory and applications, International Academic Press, USA 347 (2007).
	
	
	
	
	\bibitem{bib12} B. Davvaz and V. Leoreanu-Fotea, Krasner Hyperring Theory,
	World Scientific (2023), http://dx.doi.org/10.1142/13652.
	
	
	\bibitem{bib13} D. S. Dummit and R. M. Foote, Abstract algebra, john wile \& sons, Inc., Hoboken, NJ (2004).
	
	
	
	\bibitem{bib14}  D. Freni, A new characterization of the derived hypergroup via
	strongly regular equivalences, Communications in Algebra 30 (2002) 3977–
	3989, http://dx.doi.org/10.1081/AGB-120005830.
	
	
	\bibitem{bib15}  M. Koskas, Groupoids, semi-hypergroups and hypergroups, Journal de Mathématiques Pures et Appliquées 49 (1970) 155.
	
	
	
	\bibitem{bib16} M. Krasner, A class of hyperrings and hyperfields, International Journal of Mathematics and Mathematical Sciences 6 (1983) 307–311, http://dx.doi.org/10.1155/S0161171283000265.
	
	
	
	\bibitem{bib17} F. Marty, Sur une generalization de la notion de groupe, In 8th
	congress Math. Scandinaves (1934) 45–49.
	
	
	\bibitem{bib18}  C. G. Massouros, Free and cyclic hypermodules, Annali di Matematica Pura ed Applicata 150 (1988) 153–166, http://dx.doi.org/10.1007/BF01761468.
	
	
	\bibitem{bib19} J. R. Munkres, Topology. featured titles for topology series, Prentice Hall, Incorporated 812 (2000) 813.
	
	
	
	
	\bibitem{bib20} R. Y. Sharp, Steps in commutative algebra, Cambridge
	university press (2000), http://dx.doi.org/10.1017/CBO9780511623684.
	
	
	\bibitem{bib21} T. Vougiouklis, The fundamental relation in hyperrings. the general hyperfield, In Proc. Fourth Int. Congress on Algebraic Hyperstructures
	and Applications (AHA 1990), World Scientific (1991) 203–211.
	
	
	
	
	
	
\end{thebibliography}



\end{document}